\documentclass[12pt]{amsart}
\usepackage{amscd,amssymb,graphicx}
\author{Alice Fialowski}
\address{E\"otv\"os Lor\'and University\\
Budapest, Hungary} \email{fialowsk@cs.elte.hu}
\author{Michael Penkava}
\address{University of Wisconsin\\
Eau Claire, WI 54702-4004} \email{penkavmr@uwec.edu}

\newtheorem{thm}{Theorem}[section]

\newtheorem{lma}[thm]{Lemma}
\theoremstyle{definition}

\def \ph{\varphi}

\def \diag{\operatorname {diag}}

\def \refeq#1{equation (\ref{#1})}
\def \ra{\rightarrow}

\def \hom{\mbox{\rm Hom}}
\def \ie{\hbox{\it i.e.}}

\def \tns{\otimes}
\def \mtns{\tns\cdots\tns}

\def \k{\mbox{$\mathbb K$}}

\def \C{\mbox{$\mathbb C$}}
\def \Z{\mbox{$\mathbb Z$}}

\def\zt{\mbox{$\Z_2$}}

\def\ad{\operatorname{ad}}

\def\inv{^{-1}}
\def\d{d}

\def\im{\operatorname{Im}}

\def\dl{\delta}

\def\coder{\operatorname{Coder}}
\def\ainf{\mbox{$A_\infty$}}
\def\linf{\mbox{$L_\infty$}}
\def\and{\mbox{ \rm and }}

\def\s#1{(-1)^{#1}}

\def\pha#1#2{\ph^{#1}_{#2}}

\def\dl#1#2{\delta^{#1}_{#2}}
\def\psa#1#2{\psi^{#1}_{#2}}

\def\inv{^{-1}}

\def\dbar#1{\overline{\overline{#1}}}
\def\dl{\delta+\lambda}
\def\bdl{\bar\delta+\bar\lambda}
\def\bbdl{\overline{\overline{\delta+\lambda}}}
\def\hdlm{H_{\mu,\delta+\lambda}}

\def\hdlpm{H_{\psi,\delta+\lambda,\mu}}
\def\hrest{H_\mu(\ker(D_{\dl}))}
\def\hmu{H_\mu}
\def\hpsi{H_\psi}
\def\hdlp{H_{\psi,\delta+\lambda}}
\def\cbdl#1{\left[\overline{#1}\right]}
\def\cbdlp#1{\left[\overline{\overline{#1}}\right]}
\def\ctpsi#1{\left\{\cbdl{#1}\right\}}
\def\ctmu#1{\left\{\cbdlp{#1}\right\}}
\def\grest{\mbox{$G^{\text{rest}}$}}
\def\ggen{\mbox{$G^{\text{gen}}$}}
\def\gdiag{\mbox{$G_{M,W}$}}
\def\gdiagmd{\mbox{$G_{\delta,\mu}$}}
\def\ggenmd{\mbox{$G^{\text{gen}}_{\delta,\mu}$}}
\def\gdiagmdl{\mbox{$\gdiagmd(\lambda)$}}
\def\ggenmdl{\mbox{$G_{\delta,\mu,\lambda}$}}
\def\GL{{\mathbf{GL}}}

\begin{document}
\setlength{\multlinegap}{0pt}
\title{Extensions  of associative algebras}%

\address{}%
\email{}%

\thanks{The research of the authors was partially supported by
 OTKA grants K77757 and NK72523 and by grants
 from the University of Wisconsin-Eau Claire.}%

\dedicatory{Dedicated to Murray Gerstenhaber on his $80^\text{th}$ birthday}%
\begin{abstract}
In this paper, we give a purely cohomological interpretation of the
extension problem for associative algebras; that is the problem of
extending an associative algebra by another associative algebra. We
then give a similar interpretation of infinitesimal deformations of
extensions. In particular, we consider infinitesimal deformations of
representations of an associative algebra.
\end{abstract}
\maketitle

\section{Introduction}
Extensions of Lie and associative algebras by ideals is a classical
subject \cite{ce,mac}, which has been recast in many forms and
generalized extensively \cite{GS2,GGS}, in terms of diagrams of
algebras. Deformation theory of associative algebras is still an
active subject of research \cite{bez-gin}.

Our goal in this paper is more modest. We wish to recast the
classical ideas in the modern language of codifferentials of coalgebras
introduced
in \cite{sta4}. (A \emph{codifferential} is simply an odd coderivation whose
square is zero.)
The goal is to describe the theory of extensions of
associative algebras in a more constructive approach, because our
ultimate aim is to use the extensions as a tool to construct moduli
spaces of low dimensional algebras.

The authors have been studying moduli spaces of algebras in several
recent papers, from the point of view of algebras as codifferentials
on certain coalgebras. The modern language of codifferentials makes
it possible to express the ideas involved in extensions in a more
explicit form, which makes it easier to apply the theory in
practice. In this paper, we will illustrate how to use the
presentation of the main results by giving examples of the
construction of moduli spaces of extensions. In \cite{fp12}, we use
the ideas presented here to give a construction of the moduli space
of 3-dimensional complex associative algebras.

In some recent works, \cite{fp3,fp8,Ott-Pen}, moduli spaces of low
dimensional Lie algebras have been constructed and interpreted using
versal deformations of the algebras.  These versal deformations were
constructed by analyzing the space of coderivations of the symmetric
algebra of the underlying vector space, so giving a description of
the theory of extensions in terms of codifferentials, as we do in
this paper, makes it possible to use the computational tools we have
already developed to study the moduli spaces of algebras more
effectively.

In this paper, we give a purely cohomological interpretation of the
extension problem in terms of differentials arising from the algebra
structures. We also give a classification of infinitesimal
deformations of extensions in terms of a certain triple cohomology.
Finally, we study the problem of deformations of representations of
associative algebras, also in terms of cohomology.

The results in this paper have immediate applications to the construction
of moduli spaces of associative algebras using extensions. The authors
have been using Maple worksheets developed by one of the authors and
his students, which calculate cohomology and deformations of associative
algebras. The authors have already been using these results in conjunction
with the Maple software to construct moduli spaces, and we expect that this
software will eventually be used by others for similar calculations.

In section \ref{sec2}, we recall the definition of an extension in
terms of coderivations. In section \ref{sec3} we recall the notion
of equivalence of extensions, giving a definition of a restricted
equivalence in terms of commutative diagrams. In section \ref{sec4}
we classify infinitesimal extensions, and then in section \ref{sec5}
we classify the extensions of an algebra by a fixed bimodule
structure. In section \ref{sec6} we classify the extensions of an
associative algebra by an ideal in terms of the restricted notion of
equivalence, and then we go on to classify the extensions in terms
of a more general notion of equivalence in section \ref{sec7}. In
section \ref{sec8} we give some simple examples illustrating the
application of the classification in constructing moduli spaces of
extensions. In section \ref{sec9}, we classify infinitesimal
deformations of extensions and in section \ref{sec10} we classify
infinitesimal deformations of representations.

\section{Extensions of Associative Algebras}\label{sec2}
In this paper, we study not necessarily unital associative algebras
defined over a field $\k$, which we will assume for technical
reasons does not have characteristic 2 or 3.
We refer to an exact
sequence of associative algebras
\begin{equation}
0\ra M\ra V\ra W\ra 0.
\end{equation}
as an extension of the algebra $W$ by the algebra $M$. For
convenience, we introduce the following notation for certain
subspaces of the tensor coalgebra $T(V)=\sum_{n=0}^\infty T^n(V)$ of
$V=M\oplus W$.
\begin{align*}
T^{k,0}(M,W)=&T^k(M)\\
T^{0,l}(M,W)=&T^l(W)\\
T^{k,l}(M,W)=&M\tns T^{k-1,l}(M,W)\oplus W\tns T^{k,l-1}(M,W).
\end{align*}
In other words, $T^{k,l}(M,W)$ is the subspace of  $T^{k+1}(V)$ spanned by
tensors of $k$ elements from $M$ and $l$ elements from $W$. We also introduce
a notation for certain spaces of cochains $C(V)=\hom(T(V),V)$ on $V$.
\begin{align*}
C^k=&\hom(T^k(W),W)\\
C^{k,l}=&\hom(T^{k,l}(M,W),M).
\end{align*}
Recall that $C(V)$ is identifiable with the space $\coder(T(V)$ of
coderivations of
the tensor coalgebra $T(V)$, which means that $C(V)$ has a \zt-graded
Lie bracket. We shall sometimes refer to cochains in $C(V)$ as coderivations.

In terms of the induced bracket of
cochains, we have
\begin{align*}
[C^k,C^l]\subseteq&C^{k+l-1}\\
[C^{k,l},C^{r,s}]\subseteq&C^{k+r-1,l+s}\\
[C^k,C^{r,s}]\subseteq& C^{r,k+s-1}.
\end{align*}

The algebra structure on $V$ is determined by the following odd cochains:
\begin{align*}
&\delta\in C^2=\hom(W^2,W):& \text{the algebra structure on $W$}\\
&\psi\in C^{0,2}=\hom(W^2,M):&  \text{the "cocycle" with values in $M$}\\
&\lambda\in C^{1,1}=\hom(WM\oplus MW,M):&  \text{the ``bimodule''
structure on $M$}\\
&\mu\in C^{2,0}=\hom(M^2,M):&  \text{the algebra structure on $M$}
\end{align*}
The fact that $d$ has no terms from $\hom(WM\oplus MW\oplus M^2,W)$
reflects the fact that $M$ is an ideal in $V$. The associativity
relation on $V$ is that the odd coderivation
$$d=\delta+\lambda+\rho+\mu+\psi$$ is an odd \emph{codifferential} on $T(V)$,
which simply means that $[d,d]=0$. Now, in general, we see that
$[d,d]\in\hom(V^3,V)$. By decomposing this space and considering
which parts the brackets of the terms $\delta$, $\lambda$,
$\mu$ and $\psi$ are defined on, we obtain
\begin{align}
\label{rel1} &[\delta,\delta]=0:
\text{ The algebra structure $\delta$ on $W$ is associative.}\\
\label{rel2} &[\mu,\mu]=0:
\text{ The algebra structure $\mu$ on $M$ is associative.}\\
\label{rel3} &
[\delta,\lambda]+1/2[\lambda,\lambda]+[\mu,\psi]=0:
\text{The Maurer-Cartan equation.}\\
\label{rel5} &[\mu,\lambda]=0:\!
\text{The compatibility relation.}\\
\label {rel6} &[\delta+\lambda,\psi]=0: \text{ The cocycle condition.}
\end{align}
We also have the relations
\begin{math}
[\mu,\delta]=[\psi,\psi]=0,
\end{math}
which follow automatically, and therefore are not conditions, {\it
per se}, on the structure $d$.

When $[\mu,\psi]=0$, condition \eqref{rel3} is called the
Maurer-Cartan equation (MC equation), implying that $\delta+\lambda$
is a codifferential, and in that case, condition \eqref{rel6} is
simply the condition that $\psi$ is a cocycle with respect to this
codifferential. In general, $\psi$ is not really a cocycle, because
$\delta+\lambda$ is not a codifferential. We shall refer to
condition \eqref{rel3} as the MC equation, although this
terminology is not precisely correct.

When neither $\mu$ nor $\psi$ vanish, then in general, $\lambda$
is not a bimodule structure on $M$. However, we can still interpret the
conditions in terms of a MC
formula as follows. The sum of the two algebra
structures $\delta+\mu$ is a codifferential on $V$,
and with respect to this structure $\lambda+\psi$ satisfies an
MC formula; \ie,
\begin{equation}\label{mc-associative}
[\delta+\mu,\lambda+\psi]+\tfrac12[\lambda+\psi,\lambda+\psi]=0.
\end{equation}
Thus the combined structure $\lambda+\psi$ plays the same role
as the module structure plays in the simpler case.  The moduli space
of all extensions  of the algebra structure $\delta$ on $W$ by the
algebra  structure $\mu$ on  $W$ is given by the solutions to the MC
formula \eqref{mc-associative}. This point of view is useful if we
consider extensions where $\mu$ is assumed to be a fixed algebra
structure on $M$.

We can adopt a different point of view by noticing that
$\mu+\lambda+\psi$ satisfies an MC formula with respect to the
codifferential $\delta$; \ie, that
\begin{equation}\label{mc-second}
[\delta,\mu+\lambda+\psi]+\tfrac12[\mu+\lambda+\psi,\mu+\lambda+\psi]=0.
\end{equation}
This formulation is useful if we are interested in studying the
moduli space of all extensions of $W$ by $M$, where we don't assume
any fixed multiplication structure $\mu$ on $M$.

All these facts are well-known, for example see
\cite{gers,gers1,gers2,gers3,gers4,mac}. Our purpose in making them
explicit here is to cast the ideas in the language of codifferentials.
We summarize the main results in the theorem below.
\begin{thm} Let $\delta$ be an associative algebra structure
on $W$ and $\mu$ be an associative algebra structure on $M$.
Then $\lambda\in\hom(WM\oplus WM,M)$
 and $\psi\in\hom(W^2,M)$ determine an associative algebra structure
on $M\oplus W$ precisely when the following conditions hold.
\begin{align}
[\delta,\lambda]+\tfrac12[\lambda,\lambda]+[\mu,\psi]=0\\
[\mu,\lambda]=0\\
[\delta+\lambda,\psi]=0
\end{align}
\end{thm}
\section{Equivalence of Extensions of Associative
Algebras}\label{sec3} A (restricted) equivalence of extensions of
associative algebras is given by a commutative diagram of the form
$$
\begin{CD}
0@>>>M@>>>V@>>>W@>>>0\\
@.@|@VVfV@|\\
0@>>>M@>>>V@>>>W@>>>0\\
\end{CD}
$$
where we assume that in  the bottom row, $V$ is equipped with the
codifferential $d=\delta+\mu+\lambda+\psi$, and in the top row, it
is equipped with the codifferential
$d'=\delta+\mu+\lambda'+\psi'$, and $f$ is a morphism of
associative algebras (which is necessarily an isomorphism). The
condition that $f$ is an isomorphism of algebras is simply $d'=f^*(d)=f\inv d
f$. In order for the diagram to commute, we must have
$f(m,w)=(m+\beta(w),w)$, where $\beta\in C^{0,1}=\hom(W,M)$. We
can also express $f=\exp(\beta)$, which is convenient because then we can express
$f^*=\exp(-\ad_\beta)$, so that
\begin{equation}
f^*(d)=\exp(-\ad_\beta)(d)=d+[d,\beta]+\tfrac12[[d,\beta],\beta]+\cdots.
\end{equation}
This series is actually finite, because $[[[d,\beta],\beta],\beta]=0$. Moreover,
$[\mu,\beta]\in C^{1,1}$, while
$[\delta,\beta]$, $[\lambda,\beta]$ and $[\rho,\beta]$ all lie in $C^{0,2}$. The only nonzero term in
$[[d,\beta],\beta]$ is $[[\mu,\beta],\beta]$, which also lies in $C^{0,2}$. It follows that
$
\lambda'=\lambda+[\mu,\beta]$ and $
\psi'=\psi+[\lambda+\tfrac12[[\mu,\beta],\beta]$. Thus we have shown
\begin{thm} If $d=\delta+\mu+\lambda+\psi$ and $d'=\delta+\mu+\lambda'+\psi'$ are
two extensions of an associative algebra structure $\delta$ on $W$ by an associative algebra structure
$\mu$ on $M$, then they are equivalent (in the restricted sense) precisely when there is some $\beta\in\hom(W,M)$
such that
\begin{align}\label{equivext}
\lambda'=&\lambda+[\mu,\beta]\\
\psi'=&\psi+[[\delta+\lambda+\tfrac12[\mu,\beta],\beta]\label{equivext2}.
\end{align}
\end{thm}
We will denote the group of restricted equivalences by $\grest$. Its
elements consist of the exponentials of $\beta\in C^{0,1}$.
\section{Infinitesimal extensions and infinitesimal
equivalence}\label{sec4}
 An infinitesimal extension of $\delta$ by
$\mu$ is one of the form
\begin{equation*}
d=\delta+\mu+t(\lambda+\psi),
\end{equation*}
where $t$ is an infinitesimal parameter (\ie, $t^2=0$).

The conditions for $d$ to be an infinitesimal extension are
\begin{align}
[\delta,\lambda]+[\mu,\psi]=0\label{infr1}\\
[\mu,\lambda]=0\label{infr2}\\
[\delta,\psi]=0\label{infr3}
\end{align}

If $\alpha\in C(V)$, then denote $D_\alpha=\ad_\alpha$. When
$\alpha$ is odd and $[\alpha,\alpha]=0$, then $D^2_\alpha=0$, and
$D_\alpha$ is called a \emph{coboundary operator} on $C(V)$, and
$H_\alpha=\ker(D_\alpha)/\im(D_\alpha)$ is the \emph{cohomology}
induced by $\alpha$. The image of the cocycles from $C^{k}$ in
$H_\alpha$ will be denoted by $H^k_\alpha$, and similarly, the image
of the cocycles from $C^{k,l}$ in $H_\alpha$ will be denoted by
$H^{k,l}_\alpha$. An element $\phi$ such that $D_\alpha(\phi)=0$ is
called a \emph{$D_\alpha$-cocycle}. The bracket on $C(V)$ descends
to a bracket on $H_\alpha$, so $H_\alpha$ inherits the structure of
a Lie superalgebra. Since $\delta$, $\mu$ and $\psi$ are all
codifferentials, they determine coboundary operators. In general,
$\delta+\lambda$ is not a codifferential, so $D_{\delta+\lambda}$ is
not a coboundary operator.

Note that in the conditions for $d$
to be an infinitesimal extension, there is a certain symmetry in the
roles played by the codifferentials $\delta$ and $\mu$,
in the sense that if we interchange $\delta$ with $\mu$, and
$\psi$ with $\lambda$, then the conditions remain the same.  We have
\begin{align*}
D_\delta&:C^k\ra C^{k+1}&\quad&D_\delta:C^{k,l}\ra C^{k,l+1}\\
D_\mu&:C^k\ra 0&\quad& D_\mu:C^{k,l}\ra C^{k+1,l}.
\end{align*}
Since $[\delta,\mu]=0$, it follows that $D_\delta$ and $D_\mu$ anticommute. As a consequence,
\begin{equation*}
D_\mu:\ker(D_\delta)\ra \ker(D_\delta),
\end{equation*}
so we can define the cohomology $H_\mu(\ker\delta)$ determined by the restriction of $D_\mu$ to
$\ker(D_\delta)$. For simplicity, let us denote the cohomology class  of a
$D_\mu$-cocycle $\ph$ by $\bar\ph$. Let us consider a $D_\mu$-cocycle $\lambda$ in $C^{1,1}$. Then
the
existence of a $\psi\in C^{0,2}$ such that $[\delta,\lambda]+[\mu,\psi]=0$ and
$[\delta,\psi]=0$ is equivalent to the assertion that
$\overline{[\delta,\lambda]}=0$ in $H_\mu^{1,2}(\ker(\delta))$.

Note that even though the condition for the existence of a $\psi$
depends explicitly on $\lambda$, rather than the cohomology class
$\overline{\lambda}$, if such a $\psi$ exists for a particular
$\lambda$ in $\overline{\lambda}$, then one exists for any element
in $\overline{\lambda}$. This follows because if  $\lambda$ is
replaced by $\lambda'=\lambda+[\mu,\beta]$ and $\psi$ by
$\psi'=\psi+[\delta,\beta]$, where $\beta\in C^{0,1}$, then we
obtain a new codifferential $d'=\delta+\mu+t(\lambda'+\psi')$, which
is \emph{infinitesimally equivalent} to $d$. By infinitesimal
equivalence, we mean an equivalence determined  by an
\emph{infinitesimal automorphism} $f=\exp(t\beta)$, where $\beta\in
C^{0,1}$. (Actually, this is a restricted version of infinitesimal
equivalence. We will introduce a more general notion later.) Since
$d'=f^*(d)$, it follows that $d'$ satisfies the conditions for an
infinitesimal extension.

Now consider a fixed $D_\mu$-cocycle $\lambda$ such that $\overline{[\delta,\lambda]}=0$ in
$H_\mu^{1,2}(\ker(\delta))$, and choose some $\psi$ such  that
$[\delta,\lambda]+[\mu,\psi]=0$.  If $\psi'=\psi+\tau$ is another solution, then $[\mu,\tau]=0$
and $[\delta,\tau]=0$. Now $[\mu,\delta]=0$, so the $D_\mu$-cohomology class
$\bar\delta$ is defined. In the Lie superalgebra structure on
$\hmu$, we have
$$[\bar\alpha,\bar\beta]=\overline{[\alpha,\beta]}.$$
Since $[\bar\delta,\bar\delta]=\overline{[\delta,\delta]}=0$,
$\bar\delta$ determines a coboundary operator
$D_{\bar\delta}$ on $\hmu$. Denote the cohomology of $D_{\bar\delta}$ by $H_{\mu,\delta}$, and the
cohomology class of a $D_{\bar\delta}$-cocycle $\bar\ph$ by $[\bar\ph]$. Then $[\bar\delta,\bar\tau]=0$,
so $\tau$ determines a cohomology class $[\bar\tau]$.

On the other hand, suppose that $\bar\tau$ is
any $D_{\bar\delta}$-cocycle.
Then $[\bar\delta,\bar\tau]=0$ implies that $[\delta,\tau]$ is a $D_\mu$-coboundary. Since
$[\delta,\tau]\in C^{0,3}$, this forces $[\delta,\tau]=0$. Thus, every $D_{\bar\delta}$-cocycle $\bar\tau$ determines
an extension. In other words, $\tau\in C^{0,2}$ determines an extension precisely when $[\mu,\tau]=0$
and $[\delta,\tau]=0$.

We wish to determine when two extensions $d=\delta+\mu+t(\lambda+\psi+\tau)$ and $d'=\delta+\mu+t(\lambda+\psi+\tau')$
are infinitesimally equivalent.  First, let us suppose that $[\overline{\tau'}]=[\bar\tau]$. Then
$\overline{\tau'}=\bar\tau+[\bar\delta,\bar\alpha]$, for some $\alpha\in C^{0,1}$.  Since $\tau',\tau\in C^{0,2}$, which
contains no $D_\mu$-coboundaries, it follows that $\tau'=\tau+[\delta,\alpha]$. It is easy to see that this
implies that $d'=\exp(t\alpha)^*(d)$. Thus elements of $[\bar\tau]$ give rise to infinitesimally
equivalent extensions.  The converse is also easy to see, so the equivalence classes of infinitesimal
extensions determined by $\lambda$ are classified by the $H_{\mu,\delta}^{0,2}$ cohomology classes $[\bar\tau]$.

We summarize these results in the following theorem.
\begin{thm}\label{th3}
The infinitesimal extensions of an associative algebra structure $\delta$ on $W$ by an associative algebra structure
$\mu$ on $M$ are completely classified by the cohomology classes $\bar\lambda\in H^{1,1}_\mu$
 which satisfy the formula
$$\overline{[\delta,\lambda]}=0\in H^{1,2}_{\mu,\delta}(\ker(D_\delta))$$
together with the cohomology classes $[\bar\tau]\in H_{\mu,\delta}^{0,2}$.
\end{thm}
\section{Classification of extensions of an associative algebra by a
bimodule}\label{sec5}
 In this section we consider the special case
of an extension of $W$ by a bimodule structure $\lambda$ on $M$.
This means that $\mu=0$, so the MC formula \eqref{rel3} reduces to
the usual MC formula
\begin{equation*}
[\delta,\lambda]+\tfrac12[\lambda,\lambda]=0.
\end{equation*}

Let us relate the definition of bimodule given here with the notion of left and right module structures.
Since
\begin{equation*}
C^{1,1}=\hom(WM,M)\oplus\hom(MW,M),
\end{equation*}
we can express $\lambda=\lambda_L+\lambda_R$,
where
$\lambda_L\in\hom(WM,M)$ and $\lambda_R\in\hom(MW,M)$. Then the MC formula above is equivalent to the
three conditions on $\lambda_L$ and $\lambda_R$ below.
\begin{align*}
&[\delta,\lambda_L]+\tfrac12[\lambda_L,\lambda_L]=0& \text{$\lambda_L$ is a left-module structure.}\\
&[\delta,\lambda_R]+\tfrac12[\lambda_R,\lambda_R]=0& \text{$\lambda_R$ is a right-module structure.}\\
&[\lambda_L,\lambda_R]=0& \text{The two module structures are compatible.}
\end{align*}
Thus our definition of a bimodule structure is equivalent to the usual notion of a bimodule.

Now, $\lambda$ determines a bimodule  structure precisely when $[\mu,\psi]=0$, owing to
\eqref{rel3} in the conditions for an extension.  Thus $\bar\psi$ is well defined in $H_\mu$.
Since $[\mu,\delta+\lambda]=0$, we can define $D_{\bdl}$ on $\hmu$. Moreover $D^2_{\bdl}=0$,
so we can define its cohomology $H_{\mu,\delta+\lambda}$ on $\hmu$. Now
for a $D_\mu$-cocycle $\psi\in C^{0,2}$, $[\delta+\lambda,\psi]=0$ precisely when
$[\bdl,\bar\psi]=0$, because $[\delta+\lambda,\psi]\in C^{0,3}$, which contains no $D_\mu$-coboundaries.
Therefore, a $D_\mu$-cocycle $\psi$ determines an extension iff $\bar\psi$ is a $D_{\bdl}$-cocycle.

On the other hand, $\overline{\psi'}\in [\bar\psi]$ iff
$\psi'=\psi+[\delta+\lambda,\beta]$ for some $\beta\in C^{0,1}$. But this happens precisely in the
case when the extensions determined by $\psi$ and $\psi'$ are equivalent (in the restricted sense).
Thus we have shown the following theorem.
\begin{thm}
The extensions of $\delta$ by $\mu$ determined by a fixed bimodule structure $\lambda$ are
classified by
the cohomology classes $[\bar\psi]\in H^{0,2}_\mu$.
\end{thm}
\section{Classification of restricted equivalence classes of
extensions}\label{sec6} In this section, we assume that $\delta$ and
$\mu$ are fixed associative algebra structures on $M$ and $W$,
respectively. We want to classify the equivalence classes of
extensions under the action of the group $\grest$ of restricted
equivalences given by exponentials of maps $\beta\in\hom(W,M)$.
First, note that $\delta$ is a $D_\mu$-cocycle, and $\lambda$ must
be a $D_\mu$-cocycle by condition \eqref{rel5}, so they determine
$D_\mu$-cohomology classes $\bar\delta$ and $\bar\lambda$ in $\hmu$.
If $\lambda,\psi$ determine an extension, and
$\lambda'\in\bar\lambda$, then $\lambda',\psi'$ determine an
equivalent extension, where $\lambda'$ and $\psi'$ are given by the
formulas \eqref{equivext} and \eqref{equivext2}.  Moreover,
condition \eqref{rel3} yields the MC formula
\begin{equation}\label{mcind}
[\bar\delta,\bar\lambda]+\tfrac12[\bar\lambda,\bar\lambda]=0,
\end{equation}
which means that given a representative $\lambda$ of a cohomology class $\bar\lambda$,
there is a $\psi$ satisfying \eqref{rel3} precisely when
$\bar\lambda$ satisfies the MC-equation for
$\bar\delta$, which is a codifferential in $\hmu$.

We also need $\psi$ to satisfy condition \eqref{rel6};
\ie, $\psi\in\ker(D_{\dl})$,
which is not automatic. However, note that since $[\mu,\delta+\lambda]=0$, $D_{\dl}$ anticommutes with $D_\mu$,
which  implies that
$D_\mu$ induces a coboundary operator on $\ker(D_{\dl})$.
Because the triple bracket of any coderivation vanishes,
$[\delta+\lambda,\delta+\lambda]\in \ker(D_{\dl})$. As a consequence, we obtain
that the existence of an extension with module structure $\lambda$
is equivalent to the condition that $[\delta+\lambda,\delta+\lambda]$ is a
$D_\mu$-coboundary in the restricted complex $\ker(D_{\dl})$. In other words, there is an
extension with module structure $\lambda$ precisely when
$\overline{[\delta+\lambda,\delta+\lambda]}=0$
in the restricted cohomology $\hrest$.

Even though the complex
$\ker(D_{\dl})$ depends on $\lambda$,
the existence of an extension with module structure
$\lambda$ depends only on the $D_\mu$-cohomology
class of $\lambda$.
Thus the assertion that $\overline{[\delta+\lambda,\delta+\lambda]}=0$ in
$\hrest$ depends only on $\bar\lambda$,
and not on the choice of a representative. Of course, the $\psi$
satisfying equation \eqref{rel3} does depend on $\lambda$.
We encountered a similar situation when
analyzing infinitesimal extensions, except that there,
one had to consider only $H_\mu(\ker(\delta))$,
instead of $\hrest$.

If $\bar\lambda$ satisfies \refeq{mcind} in $\hmu$, then $D_{\bdl}^2=0$,
so we can define an associated cohomology,
which we denote by
$\hdlm$. If $\bar\ph$ is a $D_{\bdl}$-cocycle, then denote its cohomology class in $\hdlm$ by
$[\bar\ph]$. Note that \refeq{mcind} is satisfied whenever
$\overline{[\delta+\lambda,\delta+\lambda]}=0$ in
$\hrest$.

Now fix $\lambda$ and $\psi$ determining an extension.
Suppose $\lambda$ and $\psi'$ also determines an extension,
and let $\tau=\psi'-\psi$. Then it follows that $[\mu,\tau]=0$  and
$[\delta+\lambda,\tau]=0$. Thus $\bar\tau$ is a $D_{\bdl}$-cocycle.
Since $\tau\in C^{0,2}$, the condition $D_{\bdl}(\bar\tau)=0$ is equivalent
to the conditions $[\mu,\tau]=0$ and
$[\delta+\lambda,\tau]=0$. Clearly, if $D_{\bdl}(\bar\tau)=0$,
then $\psi'=\psi+\tau$ determines an extension.
Thus the set of extensions with a fixed $\lambda$ are determined
by the $D_{\bdl}$-cocycles $\bar\tau$.

We wish to determine when two extensions
$d=\delta+\mu+\lambda+\psi+\tau$ and
$d=\delta+\mu+\lambda+\psi+\tau'$ are equivalent.  If
$\overline{\tau'}\in [\bar\tau]$, then
$\overline{\tau'}=\bar\tau+[\bdl,\bar\beta]$, for some $\beta\in
C^{0,1}$, and since $\tau\in C^{0,1}$ which contains no
$D_\mu$-coboundaries, $\tau'=\bar\tau+[\delta+\lambda,\beta]$. It
follows that $d'=\exp(\beta)^*(d)$, so the extensions are
equivalent. Conversely, if $d'=\exp(\beta)^*(d)$, then
$[\mu,\beta]=0$ and  $\tau'=\tau+[\dl,\beta]$, so
$\overline{\tau'}\in[\bar\tau]$.
\begin{thm}\label{th5}
The equivalence classes of extensions of the associative algebra
structure $\delta$ on $W$ by an associative algebra structure
$\mu$ on $M$ under the action of the group $\grest$ of restricted
equivalences are completely classified by cohomology classes
$\bar\lambda\in H_\mu^{1,1}$ which satisfy the condition
$$\overline{[\delta+\lambda,\delta+\lambda]}=0\in H_{\mu}^{1,2}(\ker(D_{\delta+\lambda}))$$
together with the cohomology classes $[\bar\tau]\in H_{\mu,\delta+\lambda}^{0,2}$.
\end{thm}
\section{General Equivalence Classes of Extensions}\label{sec7}
In the standard construction of equivalence of extensions, we have assumed that the homomorphism
$f:V\ra V$ acts as the identity on $M$ and $W$.  We could consider a more general commutative
diagram of the form
$$
\begin{CD}
0@>>>M@>>>V@>>>W@>>>0\\
@.@VV\eta V@VVfV@VV\gamma V\\
0@>>>M@>>>V@>>>W@>>>0\\
\end{CD}
$$
where $\eta$ and $\gamma$ are isomorphisms.  It is easy to see
that under this circumstance, if $d'$ is the codifferential on the
top line, and $d$ is the one below, then $\eta^*(\mu)=\mu'$ and
$\gamma^*(\delta)=\delta'$.  Therefore, if one is interested in
studying the most general moduli space of all possible extensions
of all codifferentials on $M$ and $W$, where equivalence of
elements is given by diagrams above, then for two extensions to be
equivalent, $\mu'$ must be equivalent to $\mu$ as a codifferential
on $M$, and $\delta'$ must be equivalent to $\delta$ as a
codifferential on $W$, with respect to the action of the
automorphism group $\GL(M)$ on $M$ and $\GL(W)$ on $W$.

Thus, in classifying the elements of the moduli space, we first
have to consider equivalence classes of codifferentials on $M$ and
$W$. As a consequence, after making such a choice,  we need only
consider diagrams which preserve $\mu$ and $\delta$;  in other
words, we can assume  that $\eta^*(\mu)=\mu$ and  that
$\gamma^*(\delta)=\delta$.

Next note that we can always decompose a general extension diagram into one of the form
$$
\begin{CD}
0@>>>M@>>>V@>>>W@>>>0\\
@.@|@VVfV@|\\
0@>>>M@>>>V@>>>W@>>>0\\
@.@VV\eta V@VVg=(\eta,\gamma)V@VV\gamma V\\
0@>>>M@>>>V@>>>W@>>>0\\
\end{CD}
$$
where $f=\exp(\beta)$, and $g=(\eta,\gamma)$ is an element of the
group $\gdiag$ consisting of block diagonal matrices. The group
$\ggen$ of general equivalences is just the group of block upper
triangular matrices, and is the semidirect product of $\grest$
with $\gdiag$; that is, $\ggen=\gdiag\rtimes\grest$. In fact, if
$g\in\gdiag$, then $g\inv\exp(\beta)g=\exp(g^*(\beta))$.

The group $\gdiag$ acts in a simple manner on cochains. If
$g\in\gdiag$, then $g^*(C^{k,l})\subseteq C^{k,l}$ and
$g^*(C^k)\subseteq C^k$. Since $g^*D_\mu=D_{g^*(\mu)}g^*$, the
action induces a map $$g^*:H_\mu\ra H_{g^*(\mu)},$$ given by
$g^*(\bar\ph)=\overline{g^*(\ph)}$. Similarly,
$g^*D_{\dl}=D_{g^*(\delta)+g^*(\lambda)}g^*$, so we obtain a map
$$g^*:H_{\mu,\dl}\ra H_{g^*(\mu),g^*(\delta)+g^*(\lambda)},$$
given by $g([\bar\ph])=[\overline{g^*(\ph)}]$.

Let $\gdiagmd$ be the subgroup of $\gdiag$ consisting of those
elements $g$ satisfying $g^*(\mu)=\mu$ and $g^*(\delta)=\delta$.
Then $\gdiagmd$ acts on $H_\mu$, and induces a map $H_{\mu,\dl}\ra
H_{\mu,\delta+g^*(\lambda)}$. Let $\gdiagmdl$ be the subgroup of $g$
in $\gdiagmd$ such that $g^*(\lambda)=\lambda$. Thus $\gdiagmdl$
acts on both on $H_\mu$ and $H_{\mu,\dl}$.

It is easy to study the behaviour of elements in $\gdiagmd$ on
extensions. If $\lambda$ gives an extension and $g\in\gdiagmd$,
then any element $\lambda'\in g^*(\bar\lambda)$ will determine an
equivalent extension, and thus equivalence classes of
$\bar\lambda$ under the action of the group $\gdiagmd$ correspond
to equivalent extensions.

Now suppose that $\lambda$, $\psi$ gives an extension, and
$\bar\tau$ is a $D_{\bdl}$-cocycle. If $g\in\gdiagmdl$, then
$$g^*(\psi+\tau)=\psi+g^*(\psi)-\psi + g^*(\tau),$$ so that
$$[\bar\tau]\mapsto [\overline{g^*(\psi)-\psi +g^*(\tau)}]$$
determines an action of $\gdiagmdl$ on $H_{\mu,\dl}$ whose
equivalence classes determine equivalent representations.

To understand the action of $\ggen$ on extensions, first note that
any element $h\in\ggen$ can be expressed uniquely in the form
$h=g\exp(\beta)$ where $g\in\gdiag$. If $d'=h^*(d)$, for an
extension $d$, then we compute the components of the extension
$d'$ as follows.
\begin{align*}
\delta'&=g^*(\delta)\\
\mu'&=g^*(\mu)\\
\lambda'&=g^*(\lambda)+[\mu',\beta]\\
\psi'&=g^*(\psi)+[\delta'+\lambda'-\tfrac12[\mu',\beta],\beta].
\end{align*}
Clearly, $\delta'=\delta$ and $\mu'=\mu$ precisely when
$g\in\gdiagmd$.   Define the group $\ggenmd$ to be the subgroup of
$\ggen$ consisting of those $h=g\exp(\beta)$ such that
$g^*(\delta)=\delta$ and $g^*(\mu)=\mu$. In other words,
$\ggenmd=\gdiagmd\rtimes\grest$.

Define $\ggenmdl$ to be the subgroup of $\ggenmd$ consisting of
those $h$ such that $\lambda=g^*(\lambda)+[\mu,\beta]$.
$\ggenmdl$ does not have a a simple decomposition in terms of
$\gdiagmdl$, because the condition $\lambda=\lambda+[\mu,\beta]$ does not
force $g\in\gdiagmdl$. However, we can still define an action of
$\ggenmdl$ on $H_{\mu,\dl}^{0,2}$ by
$$[\bar\tau]\ra [\overline{
g^*(\psi)-\psi+g^*(\tau)+[\delta+\lambda-\tfrac12[\mu,\beta],\beta]}],$$
whose equivalence classes determine equivalent representations. Note
that for any element $\ph$ in $C^{0,2}$, $g^*(\ph)=h^*(\ph)$, so we
can use $h$ in place of $g$ in the formula above.

Note that $
0=[\mu,\tau]=g^*([\mu,\tau])=[\mu,g^*(\tau)]$,
so $\overline{g^*(\tau)}$ is well defined. Moreover,
\begin{align*}
0=[\delta+\lambda,\tau]&=g^*([\delta+\lambda,\tau])=
[\delta+g^*(\lambda),g^*(\tau)]\\&=
[\delta+\lambda-[\mu,\beta],g^*(\tau)]
=[\delta+\lambda,g^*(\tau)]-[[\mu,[\beta,g^*(\tau)]
\\&=[\delta+\lambda,g^*(\tau)]
\end{align*}
so $[\overline{g^*(\tau)}]$ is also well defined.
Thus we can define $g^*([\bar\tau])=[\overline{g^*(\tau)}]$.

In many applications, given a $\lambda$ for which a solution to the MC equation exists,
one actually has $[\mu,\psi]=0$, and thus we can choose $\psi=0$ as a solution. It also
is common that the only solutions to $\lambda=g^*(\lambda)+[\mu,\beta]$ satisfy
$[\mu,\beta]=0$. In this case, one has the much simpler formula
$[\bar\tau]\mapsto g^*([\bar\tau])$.

\begin{thm}\label{th6}
The equivalence classes of extensions of $W$ by $M$ under the
action of the group $\ggen$ are classified by the following data:
\begin{enumerate}
\item Equivalence classes of codifferentials $\delta$ on $W$ under the action
$\GL(W)$.
\item Equivalence classes of codifferentials
$\mu$ on $M$ under the action of the group $\GL(M)$.
\item Equivalence classes of $D_\mu$-cohomology classes $\bar\lambda\in H_\mu^{1,1}$
which satisfy the MC-equation
$$\overline{[\delta+\lambda,\delta+\lambda]}=0\in H^{1,2}_\mu(\ker(D_{\dl}))$$
under the action of the group $\gdiagmd$ on $\hmu$.
\item Equivalence classes of $D_{\bdl}$-cohomology classes $[\bar\tau]\in H_{\mu,\dl}^{0,2}$
under the action of the group $\ggenmdl$.
\end{enumerate}
\end{thm}
We are more interested in the moduli space of extensions of $W$ by $M$ preserving fixed
codifferentials on these spaces.
\begin{thm}\label{th7}
The equivalence classes of extensions of a codifferential $\delta$
on $W$ by a codifferential $\mu$ on $M$ under the action of the
group $\ggenmd$ are classified by the following data:
\begin{enumerate}
\item Equivalence classes of $D_\mu$-cohomology classes $\bar\lambda\in H_\mu^{1,1}$
which satisfy the MC-equation
$$\overline{[\delta+\lambda,\delta+\lambda]}=0\in H^{1,2}_\mu(\ker(D_{\dl}))$$
under the action of the group $\gdiagmd$ on $\hmu$.
\item Equivalence classes of $D_{\bdl}$-cohomology classes $[\bar\tau]\in H_{\mu,\dl}^{0,2}$
under the action of the group $\ggenmdl$.
\end{enumerate}
\end{thm}
To illustrate why this more general notion of equivalence is useful, we give some simple examples
of extensions of associative algebras. For simplicity, we assume that the base field is $\C$ in
all our examples.
\section{Simple examples of extensions of associative
algebras}\label{sec8} The notion of a Lie superalgebra is
expressable in terms of coderivations on the symmetric coalgebra of
a \zt-graded vector space. These superalgebras have been well
studied.  The corresponding notion of an associative algebra on a
\zt-graded vector space is not as well known.  One reason for this
might be that the definition of an associative algebra on a graded
vector space is the same as for a non-graded space; the
associativity relation does not pick up any signs as happens with
the graded Jacobi identity on a superspace. At first glance, it does
not appear that there is any reason to study such super associative
algebras.  However, the notion of an \ainf\ algebra, which
generalizes the idea of an associative algebra, naturally arises in
the \zt-graded setting. The study of associative algebra structures
on \zt-graded spaces is really the first step in the study of \ainf\
algebras.

The manner in which the \zt-grading appears in the classification of
associative algebra structures on a \zt-graded space is in terms of
the parity of the multiplication, which is always required to be
even, and in the parity of automorphisms of the vector space, which
are also required to be even maps.  Thus, for a \zt-graded space,
not every multiplication on the underlying ungraded space is
allowed, and only certain automorphisms of the space are allowed.
This effects  the moduli space in two ways. First, there are fewer
codifferentials, and secondly, the set of equivalences is
restricted, so it is not obvious whether the moduli space of
\zt-graded associative algebras is larger or smaller than the moduli
space on the associated ungraded space. In fact, there is  a natural
map between the moduli space of \zt-graded associative algebras and
the moduli space of ungraded algebras. In general, this map may be
neither surjective, nor injective.

Because it is convenient to work in the parity reversed model,
an ungraded vector space will correspond
to a completely odd space in our setting.
Moreover, the associativity relation picks up signs in the
parity reversed model.
In fact, if $d$ is an odd codifferential in $C^2(W)=\hom(T^2(W),W)$, then
the associativity relation becomes
\begin{equation}
d(d(a,b),c)+\s{a}d(a,d(b,c))
\end{equation}
Note that when $V$ is completely odd, then $\s{a}=-1$ for all $a$,
which gives the usual associativity relation.  Nevertheless, the
relation above gives the usual associativity relation on the parity
reversion $V=\Pi(W)$, because the induced multiplication is given by
$$m(x,y)=\s{x}\pi d(\pi\inv(x),\pi\inv(y)).$$

If $V=\langle e_1,\ldots, e_m\rangle$ is a \zt-graded space, then a
basis for the $n$-cochain space $C^n(V)=\hom(T^n(V),V)$ is given by
the coderivations $\ph^I_i$, where $I=(i_1,\ldots,i_n)$ is a
multi-index, and
\begin{equation*}
\ph^I_i(e_J)=\delta^I_Je_i.
\end{equation*}
Here $e_J=e_{j_1}\mtns e_{j_n}$ is a basis element of $T^n(V)$,
determined by the multi-index $J$
and $\delta^I_J$ is the Kronecker delta. When $\ph^I_i$ is odd, we
denote it by $\psi^I_i$ to emphasize this fact.

If $V=\langle e_1\rangle$ is a 1-dimensional odd vector space, (we
denote its dimension by $0|1$) then it has only one nontrivial odd
codifferential of degree 2, namely $d=\psi^{1,1}_1$. In other words
$d(e_1,e_1)=e_1$. On the other hand, an even 1-dimensional vector space
($1|0$-dimensional) has no nontrivial odd codifferentials. As the first
case of examples of the theory of extensions, we will study extensions
of 1-dimensional spaces by 1-dimensional spaces, as the construction is
very easy.

\subsection{Extensions where $\dim(V)=0|2$}
The classification of 2-di\-men\-sional asso\-ciative algebras dates
back at least to \cite{pie}. These algebras form the moduli space of
$0|2$-dimensional associative algebras. Let $V=\langle f_1,f_2\rangle$,
where $f_1$ and $f_2$ are both odd.
\begin{align*}
d_1&=\psa{1,1}1+\psa{2,2}2\\
d_2&=\psa{1,1}1+\psa{1,2}2\\
d_3&=\psa{1,1}1+\psa{2,1}2\\
d_4&=\psa{1,1}1+\psa{1,2}2+\psa{2,1}2\\
d_5&=\psa{1,1}1\\
d_6&=\psa{1,1}2
\end{align*}
All of these codifferentials arises as an extension of a
$0|1$-dimensional space $W$ by a $0|1$-dimensional space $M$,
because there are no simple, complex 2-dimensional complex
associative algebras. We will express $M=\langle f_2\rangle$ and
$W=\langle f_1\rangle$. Then we can have $\delta=\psa{1,1}1$ or
$\delta=0$ up to equivalence, and similarly $\mu=\psa{2,2}2$ or
$\mu=0$ up to equivalence. We can express
$\lambda=a\psa{1,2}2+b\psa{2,1}2$, $\psi=c\psa{1,1}2$ and
$\beta=x\pha12$. Then we have
\begin{align*}
\tfrac12[\lambda,\lambda]&=-a^2\psa{1,1,2}2+b^2\psa{2,1,1}2\\
[\lambda,\psi]&=(b-a)c\psa{1,1,1}2.
\end{align*}
These formula will prove useful in calculating the extensions.
\subsubsection{Case 1: $\delta=\mu=0$}
Since $\mu=0$, the module-algebra compatibility relation \eqref{rel5} is automatic.
The MC formula reduces to $[\lambda,\lambda]=0$, which forces
$\lambda=0$. The group $G_{\mu,\delta}$ consists of the diagonal
automorphisms $g=\diag(r,u)$, where $ru\ne0$, and thus the group
$\ggenmdl$ consists of products of an arbitrary diagonal matrix and
an exponential of the form $\exp(x\pha12)$. We can choose $\psi=0$,
and $\tau=c\psa{1,1}2$. Then $g^*(\tau)=\tfrac u{r^2}\tau$, so when
$c\ne0$, we can assume that $c=1$. Note that $\beta$ plays no role
in determining the equivalence class of $\tau$, because
$[\delta+\lambda+\tfrac12[\mu,\beta],\beta]=0$. Therefore, the
moduli space consists of the single element $d=\psa{1,1}2$, which is
the codifferential $d_6$ in the moduli space on $V$.
\subsubsection{Case 2: $\delta=\psa{1,1}1$ and $\mu=0$}
The group $\gdiagmd$ consists of diagonal matrices of the form
$g=\diag(1,u)$. Since $[\delta,\psi]=0$, the cocycle condition
\eqref{rel6} forces $\psi=0$, unless $a=b$.  The MC formula reduces to
\begin{equation*}
0=[\delta,\lambda]+\tfrac12[\lambda,\lambda]=(a-a^2)\psa{1,1,2}2-(b-b^2)\psa{2,1,1}2,
\end{equation*}
so $a$ and $b$ can either be 0 or 1, leading to 4 subcases.

\indent{\it Subcase 1: $\lambda=0$.}\\
We can choose $\psi=0$, and thus $\tau=c\psa{1,1}2$. If
$g=\diag(1,u)$, then $g^*(\tau)=\tfrac1u\tau$, so
\begin{equation*}
g^*(\psi)-\psi+g^*(\tau)+[\delta+\lambda-\tfrac12[\mu,\beta],\beta]=(\tfrac
cu-x)\psa{1,1}2.
\end{equation*}
As a consequence, the only case that we have to consider is
$\tau=0$, so we obtain the extension given by the codifferential
$d_5$ on the list of codifferentials on $V$.

\indent{\it Subcase 2: $\lambda=\psa{1,2}2+\psa{2,1}2$.}\\
Note that $g^*(\lambda)=\lambda$ for any $g\in\gdiagmd$. From this
we deduce that $\ggenmdl=\gdiagmd\rtimes\grest$.  As in the previous case,
we can choose $\psi=0$ and we have
$[\delta+\lambda,\beta]=\psa{1,1}2$, so $\tau=0$ gives the only
equivalence class under the action of $\ggenmdl$ on
$H_{\mu,\dl}^{0,2}$. Thus we obtain only 1 codifferential, which is
$d_4$ on the list.

\indent{\it Subcase 3: $\lambda=\psa{1,2}2$.}\\ We must have
$\psi=0$, and we we obtain the codifferential $d_2$.

\indent{\it Subcase 4: $\lambda=\psa{2,1}2$.}\\ Again, we have
$\psi=0$, and we obtain the codifferential $d_3$.

Thus we obtain exactly the 4 codifferentials $d_2$, $d_3$, $d_4$ and
$d_5$ as extensions of $\delta=\psa{1,1}1$ by $\mu=0$.
\subsubsection{Case 3: $\delta=0$ and $\mu=\psa{2,2}2$.}
Since $[\mu,\lambda]=(b-a)\psa{2,1,2}2$ must vanish by condition
\refeq{rel5}, this forces $a=b$ and
$\lambda=a\psa{1,2}2+a\psa{2,1}2$. But
$[\mu,\beta]=c\psa{1,2}2+c\psa{2,1}2$ which means $\lambda$ is a
$D_\mu$-coboundary. Thus we can assume $\lambda=0$.

Moreover the MC formula forces $\psi=0$, since
$[\mu,\psi]=c(\psa{1,1,2}2-\psa{2,1,1}2)$. Thus we obtain the
codifferential $d=\psa{2,2}2$.  This codifferential does not appear
on the list of codifferentials on $V$, but is equivalent to $d_5$
under an obvious change of basis.
\subsubsection{Case 4:$\delta=\psa{1,1}1$ and $\mu=\psa{2,2}2$.}
By the same argument as the previous case, we must assume that
$\lambda=0$ and this forces $\psi=0$. Thus we obtain the
codifferential $d=\psa{1,1}1+\psa{2,2}2$, which is $d_1$ on the list
of codifferentials on $V$.

Thus  every codifferential on the $0|2$-dimensional vector space arises as an extension.
This fact is not surprising, since we know that no 2-dimensional associative algebra is simple,
so there must be a nontrivial ideal. In fact, the construction of the moduli space
can be made much simpler if one takes into account the Fundamental Theorem of
Finite Dimensional Algebras, which says that any nonnilpotent algebra is a semidirect
sum of a nilpotent ideal and a semisimple algebra. Our motive here was to illustrate
the ideas, rather than to give the simplest construction.

\subsection{Extensions where $\dim V=1|1$}
If $V=\langle e_1,e_2\rangle$, where $e_1$ is even and $e_2$ is odd
(the convention is to list the even elements in a basis first), then
the moduli space of associative algebras on $V$ contains exactly 6
nonequivalent codifferentials, just as in the $0|2$-dimensional case:
\begin{align*}
d_1&=\psa{2,2}2-\psa{1,1}2-\psa{1,2}1+\psa{2,1}1\\
d_2&=\psa{2,2}2-\psa{1,2}1\\
d_3&=\psa{2,2}2+\psa{2,1}1\\
d_4&=\psa{2,2}2-\psa{1,2}1+\psa{2,1}1\\
d_5&=\psa{2,2}2\\
d_6&=\psa{1,1}2
\end{align*}
The first element in the list is a simple $1|1$-dimensional algebra, so is not
obtainable as an extension.
We will show that the other five of these codifferentials
arise either as an extension of a 1-dimensional
odd algebra by a 1-dimensional even one or vice-versa.
\subsubsection{Extension of a $1|0$-space by a $0|1$ space}
Let $V=\langle e_1,e_2\rangle$, where $W=\langle e_2\rangle$ is odd and
$M=\langle e_1\rangle$ is an even
vector space.  Then the only nontrivial multiplication on $W$ (up to equivalence) is
$\delta=\psa{2,2}2$.  Because $M$ is even, we must have $\mu=0$.
Since $C^{0,2}=\langle \pha{2,2}1\rangle$
is even, the cocycle $\psi$ must vanish.  On the other hand, we have
$C^{1,1}=\langle \psa{1,2}1,\psa{2,1}1\rangle$ is completely odd, so the module structure
$\lambda=a\psa{1,2}1+b\psa{2,1}1$.
Then
\begin{align*}
[\delta,\lambda]&=a\psa{1,2,2}1+b\psa{2,2,1}1\\
\tfrac12[\lambda,\lambda]&=a^2\psa{1,2,2}1-b^2\psa{2,2,1}1
\end{align*}
Thus, to satisfy the MC formula \eqref{rel3}, we must have $a+a^2=0$ and $b-b^2=0$.
This gives 4 solutions,
corresponding to the codifferentials $d_2$-$d_5$ in the
list of codifferentials on a $1|1$-dimensional
space.  Note that there are no nontrivial restricted automorphisms,
because $W$ and $M$ have opposite
parity, so any map $\beta\in C^{0,1}$ must be odd.

If $\delta=0$ is chosen instead, then it follows that
$\lambda=0$ as well, so we don't obtain any nontrivial
extensions.

The group $G(\mu,\delta)$ consists of maps of the
form $g=(c,1)$, where $c$ is a nonzero constant,
because only the identity map of $W$ preserves $\delta$.\
Any such $g$ acts trivially on $\lambda$, that is $g^*(\lambda)=\lambda$.

\subsubsection{Extensions of a $0|1$-space by a $1|0$ space}
Let $V=\langle e_1, e_2\rangle$ where $W=\langle e_1\rangle$ and $M=\langle e_2\rangle$.
This time $\delta=0$,
and $\mu=\psa{2,2}2$ gives the nontrivial equivalence class of structures on $M$.
We must have
$\lambda=0$ because $C^{1,1}$ is even, but we can have $\psi=a\psa{1,1}2$.
However, since
$[\psa{1,1}2,\psa{2,2}2]=\psa{1,1,2}2-\psa{2,1,1}2$, the MC equation forces $\psi=0$. Thus we
only obtain the codifferential $d=\psa{2,2}2$ as an extension when $\mu$ is nontrivial.
However,
when $\mu=0$, $\psi=a\psa{1,1}2$ need not vanish.
In this case, the group $G_{\mu,\delta}$ is given by
automorphisms of the form $g=(b,c)$ where $b$, and $c$ are arbitrary nonzero constants.
Then
$g^*(\psi)=b^{-2}\psi$,
we can assume that $a=1$ and we obtain the codifferential $d_6$.
Already in this example, we see the
importance of using the general notion of equivalence,
because with only the notion of restricted
equivalence, the codifferentials $a\psa{1,1}2$ are nonequivalent,
and thus would be considered as
distinct extensions.

\section{Infinitesimal deformations of extensions of associative algebras}\label{sec9}
A natural question that arises when studying the moduli spaces arising from extensions is how to fit the
moduli together as a space, and to answer that question, one needs to have a notion of how to move around
in the moduli space.  This notion is precisely the idea of deformations, in this case, deformations of
the extensions.  We will classify the infinitesimal deformations of an extension.

Let $d=\delta+\mu+\lambda+\psi$ be an extension, and consider the
infinitesimal deformation
$$
d_t=d+t(\eta+\zeta)
$$
of this extension, where $\eta\in C^{1,1}$ represents a deformation
of the $\lambda$ structure, and $\zeta\in C^{0,2}$ gives a
deformation of the $\psi$ structure. In this section, we don't
consider deformations which involve deforming the $\delta$ or $\mu$
structure. The infinitesimal condition  is that $t^2=0$, in which
case, as usual, the condition for $\eta,\zeta$ to determine a
deformation is, infinitesimally, that $[d,\eta+\zeta]=0$. We split
this one condition up into the four conditions below.
\begin{align}
&[\delta+\lambda,\eta]+[\mu,\zeta]=0\label{cond1}\\
&[\delta+\lambda,\zeta]+[\psi,\eta]=0\label{cond2}\\
&[\mu,\eta]=0\label{cond3}\\
&[\psi,\zeta]=0\label{cond4}
\end{align}
These conditions are symmetric in the roles of $\psi$ and $\mu$, but
this symmetry is a bit misleading. For example the condition \eqref{cond4} is automatic for
$\zeta\in C^{0,2}$, but condition \eqref{cond3} is not automatic for $\eta\in C^{1,1}$.

Note that since $\eta$ is a $D_\mu$-cocycle, $\bar\eta$ is well defined, and condition \eqref{cond1} implies
that $\bar\eta$ is a $D_{\bdl}$-cocycle.  Since $[\psi,\psi]=0$, it determines a coboundary
operator $D_\psi$ as well.  Denote the $D_\psi$-cohomology class of a $D_\psi$-cocycle $\ph$ by
$\dbar\ph$ and the set of cohomology classes by
$\hpsi$. Note that $\hpsi$ inherits the structure of a Lie algebra.

Since $[\delta+\lambda,\psi]=0$, it follows that $\bbdl$ is well defined.
Moreover, we have
\begin{align*}
[\dbar{\delta+\lambda},[\dbar{\delta+\lambda},\dbar\ph]]&=\dbar{[\delta+\lambda,[\delta+\lambda,\ph]]}=
\dbar{[\tfrac12[\delta+\lambda,\delta+\lambda],\ph]}
\\&=-\dbar{[[\mu,\psi],\ph]}=-\dbar{[\mu,[\psi,\ph]]}=0,
\end{align*}
so $D_{\bbdl}$ is a differential on $\hpsi$.
Denote the cohomology class of a $D_{\bbdl}$-cocycle $\dbar\ph$ by
$\cbdlp\ph$ and the set of cohomology classes by $\hdlp$. Note that $\hdlp$ inherits the structure
of a Lie algebra.

We first remark that conditions \eqref{cond1} and \eqref{cond3} imply that
$\left[\bar\eta\right]$ is well defined, and \eqref{cond2} and \eqref{cond4} imply that
$\cbdlp\zeta$ is well defined.

Next we introduce an action of $D_\psi$ on $\hdlm$.
It is not possible to extend the operation of bracketing with $\psi$ to the $D_\mu$-cohomology,
because $[\mu,\psi]\ne0$.  Moreover, even if $[\mu,\ph]=0$, it does not follow that $[\mu,[\psi,\ph]]=0$.
However, we can  extend  the bracket to $\hdlm$ as follows.
A cohomology class $\left[\bar\ph\right]$, is given by a $\ph$ such that $[\mu,\ph]=0$ and
$[\delta+\lambda,\ph]=[\mu,\beta]$ for some $\beta$. Note that
\begin{align*}
[\mu,[\psi,\ph]]&=[[\mu,\psi],\ph]=-[\delta+\lambda,[\delta+\lambda,\ph]]
\\&=-[\delta+\lambda,[\mu,\beta]]=[\mu,[\delta+\lambda,\beta]].
\end{align*}
In \cite{fp9}, it was shown that the operator $D_{\psi}$ on $H_{\mu,\dl}$ given by
\begin{equation*}
D_\psi([\bar\ph])=\left[\overline{[\psi,\ph]-[\delta+\lambda,\beta]}\right],
\end{equation*}
where $\beta$ is any solution to $[\delta+\lambda,\ph]=[\mu,\beta]$, is well defined,
and that $D^2_\psi=0$. Moreover, if $H_{\mu,\dl,\psi}$ denotes the associated cohomology,
then the bracket on $H_{\mu,\dl}$ descends to a bracket on $H_{\mu,\dl,\psi}$, equipping it
with the structure of a Lie superalgebra. Let us denote the
the $D_\psi$-cohomology class of a $D_\psi$-cocycle  $\cbdl{\ph}$ by $\ctpsi{\ph}$.

In a very similar manner, one can show that one can define $D_\mu$  on
$\hdlp$ by
$$
D_\mu(\cbdlp\ph)=\cbdlp{[\mu,\ph]-[\delta+\lambda,\beta]},
$$
where $\beta$ is any coderivation satisfying $[\delta+\lambda,\ph]=[\psi,\beta]$. Then $D_\mu$ is a Lie
algebra morphism on $\hdlp$ whose square is zero, and we denote the resulting  cohomology by $\hdlpm$
and the cohomology class of a $D_\mu$-cocycle $\cbdlp{\ph}$ by $\ctmu{\ph}$.

We will call $H_{\mu,\delta+\lambda,\psi}$ and
$H_{\psi,\delta+\lambda,\mu}$ \emph{triple} cohomology groups. It
turns out that the first one will play a more important role in the
classification of infinitesimal deformations of extensions. The
following lemma was proved in \cite{fp9}.
\begin{lma}\label{eta-zeta}
Let $d=\delta+\mu+\lambda+\psi$ be an extension of the
codifferentials $\delta$ on $W$ by $\mu$ on $M$, that
$\eta\in\hom(MW,M)$ and $\zeta\in\hom(W^2,M)$. If
$$d_t=d+t(\eta+\zeta)$$
determines an infinitesimal deformation of $d$ then
\begin{enumerate}
\item $\ctpsi\eta$ is well defined.
\item $\ctmu\zeta$ is well defined.
\end{enumerate}
\end{lma}
It turns out that infinitesimal deformations can be characterized in
terms of the triple cohomology $H_{\mu,\dl,\psi}$ alone. The
following theorem, proved in \cite{fp9}, gives a condition for an
infinitesimal deformation to exist, depending on $\eta$ alone, and
classifies all such deformations.

\begin{thm}\label{th8}
Let $d=\delta+\mu+\lambda+\psi$ be an extension of the
codifferentials $\delta$ on $W$ by $\mu$ on $M$.

An element $\eta\in C^{1,1}$ gives rise to an
infinitesimal deformation for some $\zeta\in C^{0,2}$
if and only if the triple cohomology class $\{[\bar\eta]\}$ in $H_{\mu,\dl,\psi}^{1,1}$ is well defined.
In this case, if $\zeta\in C^{0,2}$ is
any coderivation such that $\eta$, $\zeta$ determine an infinitesimal deformation, then
$\zeta'=\zeta+\tau$ determines another infinitesimal deformation if and only if the double cohomology
class $[\bar\tau]$ is well defined in $H_{\mu,dl,\psi}^{0,2}$.

Moreover the infinitesimal equivalence classes of infinitesimal deformations are classified by
the triple cohomology
classes $\{[\bar\eta]\}\in H_{\mu,\dl,\psi}^{0,1}$ and $\{[\bar\tau]\}\in H_{\mu,\dl,\psi}^{0,2}$.
\end{thm}

\section{Infinitesimal Deformations of Representations}\label{sec10}
In this section, we give a complete classification of infinitesimal
deformations of representations of associative algebras

Let $M$ be an associative algebra with multiplication $\mu$, which
is also a module over $W$. In other words, we are studying an
extension of $W$ by $M$ for which the cocycle $\psi$ vanishes. There
are two interesting problems we could study.
\begin{enumerate}
\item Allow the module structure $\lambda$ and the algebra structure $\delta$ to vary,
but keep $\mu$ fixed. This case includes the study of deformations of a module structure
where the module does not have an algebra structure.
\item Allow the module structure $\lambda$ and the multiplication $\mu$ to vary,
but keep the algebra structure $\delta$ fixed.
\end{enumerate}
In both of these scenarios, we think of the structures on $M$ and
$W$ as being distinct, with interaction only through $\lambda$, so
when considering automorphisms of the structures, it is reasonable
to restrict to automorphisms of $V$ which do not mix the $W$ and $M$
terms, in other words we allow only elements of $\gdiag$.

Then we have the following maps:
\begin{align*}
&D_{\delta}:C^n\ra C^{n+1}\\
&D_{\lambda}:C^n\ra C^{1,n}\\
&D_{\delta+\lambda}:C^{k,l}\ra C^{k,l+1}\\
&D_\mu:C^{k,l}\ra C^{k+1,l}.
\end{align*}
In the setup of this problem, we only are interested in $C^{k,l}$ for $k\ge 1$, so
we shall restrict our space of cochains in this manner.  Because of this restriction,
we note that an element in $C^{1,1}$ can be a $D_\mu$-cocycle, but never a $D_{\mu}$-
coboundary. Moreover $C^n\subseteq\ker(D_\mu)$, so an element in $C^2$ is always a
$D_\mu$-cocycle, and never a $D_\mu$-coboundary.

Because $\psi=0$, the MC-equation $[\delta,\lambda]+\tfrac12[\lambda\lambda]=0$ is
satisfied, so that $D_{\delta+\lambda}^2=0$.
Since
\begin{align*}
(D_\lambda D_\delta+D_{\delta+\lambda}D_\lambda)(\ph)=&
[\lambda,[\delta,\ph]]+[\delta+\lambda,[\lambda,\ph]]\\
=&[[\lambda,\delta],\ph]-[\delta,[\lambda,\ph]]+[\delta,[\lambda,\ph]]+[\lambda,[\lambda,\ph]]
\\=&[[\delta,\lambda],\ph]+[\tfrac12[\lambda,\lambda],\ph]=0.
\end{align*}
 we have
\begin{equation*}
D_\lambda D_\delta+D_{\delta+\lambda}D_\lambda=0.
\end{equation*}
If we denote the $D_\mu$-cohomology class of a $D_\mu$-cocycle $\ph$ by $\bar\ph$ as
usual, then since $\bar\lambda$ and $\bar\delta$ are defined, we get the following
version of this equation, applicable to the cohomology space $H_\mu$.\begin{equation*}
D_{\bar\lambda} D_{\bar\delta}+D_{\bdl}D_{\bar\lambda}=0.
\end{equation*}
As usual, let us denote the $D_{\bdl}$-cohomology class of a $D_{\bdl}$-cocycle $\bar\ph$
by $[\bar\ph]$.

Let us study the first scenario, where we allow $\lambda$ and $\delta$ to vary, in other
words, we consider
\begin{equation*}
d_t=d+t(\delta_1+\lambda_1),
\end{equation*}
where $\delta_1\in C^2$ and $\lambda_1\in C^{1,1}$ represent the variations in $\delta$ and
$\lambda$.  The infinitesimal condition $[d_t,d_t]=0$ is equivalent to the three
conditions for a deformation of a module structure:
\begin{align*}
[\delta,\delta_1]&=0\\
[\lambda,\delta_1]+[\delta+\lambda,\lambda_1]&=0\\
[\mu,\lambda_1]&=0.
\end{align*}
By the third condition above $\bar\lambda_1$ is well defined, and $\bar\delta_1$ is defined.
We claim that if $D_{\bar\delta}(\bar\delta_1)=0$, which is the first condition, then
the $D_{\bdl}$-cohomology class $[D_{\bar\lambda}(\delta_1)]$ is well defined and depends
only on the $D_\delta$-cohomology class of $\delta_1$. It is well defined because
\begin{equation*}
D_{\bdl}D_{\bar\lambda}(\bar\delta_1)=-D_{\bar\lambda}D_{\bar\delta}(\bar\delta_1)=0.
\end{equation*}
To see that it depends only on the $D_\delta$-cohomology
class of $\bar\delta$, we apply $D_{\bar\lambda}$ to a $D_{\bar\delta}$-coboundary
$D_{\bar\lambda}(\bar\ph)$ to obtain
\begin{equation*}
D_{\bar\lambda}D_{\bar\delta}(\bar\ph)=D_{\bdl}D_{\bar\lambda}(-\bar\ph),
\end{equation*}
which is a $D_{\bdl}$-coboundary.
The second condition for a deformation of the module structure implies that
$[D_{\bar\lambda}(\bar\delta)]=0$. Moreover, if this statement holds, then there is
some $\lambda_1$ such that $\delta_1$ and $\lambda_1$ determine a deformation of
the module structure.  We see that $\lambda_1'=\lambda+\tau$ is another solution precisely when
$\bar\tau$ exists and $D_{\bdl}(\bar\tau)=0$. Thus, given one solution
$\lambda_1$, the set of solutions is
determined by the $D_{\bdl}$-cocycles $\tau\in C^{1,1}$.

Now let us consider infinitesimal equivalence.  We suppose that $\alpha\in C^{1,0}$ and
$\gamma\in C^1$, and $g=\exp(t(\alpha+\beta)$. If $d_t'=g^*(\d_t)$ is given by
the cochains $\delta_1'$ and $\lambda_1'$, then we have
\begin{align*}
\delta_1'&=\delta_1+D_\delta(\gamma)\\
\lambda_1'&=\lambda_1+D_{\lambda}(\alpha+\gamma)\\
D_{\mu}(\alpha+\gamma)&=0.
\end{align*}
It follows that the set of equivalence classes of deformations are determined by
$D_\delta$ cohomology classes of $\delta_1\in C^2$. If we fix $\delta_1$ such that
$D_\delta(\delta_1)=0$ and $\lambda_1$ satisfying the rest of the conditions of a
deformation, then expressing $\tau'=\tau+D_{\lambda}(\alpha+\gamma)$.
But, since $D_\delta(\alpha+\gamma)=0$, this means
we can express $\tau'=\tau+D_{\delta+\lambda}(\alpha+\gamma)$ and $D_\mu(\alpha+\gamma)$,
which means that $\bar\tau'=\bar\tau+D_{\bdl}(\bar\alpha+\bar\gamma)$, and the solutions
for $\tau$ are given by $D_{\bdl}$-cohomology classes of
$D_\mu$-cocycles $\tau\in C^{1,1}$.
Thus we obtain
\begin{thm}\label{th9}
The infinitesimal deformations of a module $M$ with multiplication
$\mu$ over an associative algebra $\delta$, allowing the algebra
structure $\delta$ on $W$ and module structure $\lambda$ to vary are
classified by
\begin{enumerate}
\item $D_{\delta}$-cohomology classes of $D_\delta$-cocycles
$\delta_1\in C^2$ satisfying the condition
$$[D_{\bar\lambda}(\bar\delta_1)]=0$$.
\item $D_{\bdl}$-cohomology classes $[\bar\tau]$ of $D_{\bdl}$-cocycles $\bar\tau$ of
$D_\mu$-cocycles $\tau\in C^{1,1}$.
\end{enumerate}
\end{thm}

Finally, let us study the second scenario, where we allow $\lambda$ and $\mu$,
but not $\delta$, to vary.  We write $d_t=d+t(\lambda_1+\mu_1)$, where
$\lambda_1\in C^{1,1}$ is the variation of $\lambda$ and $\mu_1\in C^{2,0}$ is
the variation in $\mu$.  The Jacobi identity $[d_t,d_t]=0$ gives three conditions
for a deformation of the module structure.
\begin{align*}
&D_\mu(\mu_1)=0\\
&D_{\delta+\lambda}(\mu_1)+D_\mu(\lambda_1)=0\\
&D_{\delta+\lambda}(\lambda_1)=0.
\end{align*}
Recall that $D_\mu$ maps $\ker(D_{\delta+\lambda})$ to itself, so
$H_\mu(\ker(\delta+\lambda))$ is well defined. The first condition on a deformation
says that $\bar\mu_1$ is well defined.  We claim that in that case,
$\overline{D_{\delta+\lambda}(\mu_1)}$ is
a well defined element of $H_\mu(\ker(\delta+\lambda)$ which depends only on $\bar\mu_1$.
This is clear, since $D_{\delta+\lambda}(\mu_1)\in\ker(D_{\delta+\lambda})$,
and  $D_\mu D_{\delta+\lambda}(\mu_1)=-
D_{\delta+\lambda}D_u(\mu_1)=0$. The second condition on a deformation says simply that
$\overline{D_{\delta+\lambda}(\mu_1)}=0$, and the fact that this statement is true
in $H_\mu(\ker(\delta+\lambda))$ is the third condition.  Therefore, assuming that
$\overline{D_{\delta+\lambda}(\mu_1)}=0$, we can find a $\lambda_1$ so that all
of the conditions for a deformation are satisfied.

If $\lambda_1+\tau$ gives another solution, then $D_\mu(\tau)=0$
and $D_{\delta+\lambda}(\tau)=0$.
Because we do not allow elements of $C^{0,1}$ as cochains, these two equalities are
equivalent to  $[\bar\tau]$ being well defined.

If $d_t'=\exp(t(\alpha+\gamma))$, then we obtain the following.
\begin{align*}
[\delta,\alpha+\gamma]&=0\\
\lambda_1'&=\lambda+[\lambda,\alpha+\gamma]\\
\mu_1'&=\mu_1+[\mu,\alpha+\gamma].
\end{align*}

Thus, up to equivalence, a deformation is given by a
$D_\mu$-cohomology class $\bar\mu_1$. If we fix $\mu_1$, and then
look at the variation in $\tau$, one also obtains that up to
equivalence, the deformation is determined by the
$D_{\bdl}$-cohomology class $[\bar\tau]$ of the $D_\mu$-cohomology
class $\bar\tau$. Thus we have shown
\begin{thm}\label{th10}
The infinitesimal deformations of a module $M$ with associative
algebra structure $\mu$ over an associative algebra $\delta$
allowing the algebra structure $\mu$ on $W$ and module structure
$\lambda$ to vary are classified by
\begin{enumerate}
\item $D_\mu$ cohomology classes $\bar\mu_1$ of $D_\mu$-cocycles $\mu_1$ lying in
$C^{2,0}$.
\item $D_{\bdl}$-cohomology classes $[\bar\tau]$ of $D_{\bdl}$-cocycles $\bar\tau$ of $D_\mu$-cocycles
$\tau$ lying in $C^{1,1}$.
\end{enumerate}
\end{thm}

\bibliographystyle{amsplain}
\providecommand{\bysame}{\leavevmode\hbox to3em{\hrulefill}\thinspace}
\providecommand{\MR}{\relax\ifhmode\unskip\space\fi MR }
\providecommand{\MRhref}[2]{%
  \href{http://www.ams.org/mathscinet-getitem?mr=#1}{#2}
}
\providecommand{\href}[2]{#2}

\end{document}